\documentclass[11pt]{article}
\usepackage{amsmath,amssymb,amsthm,color}

\usepackage{comment}
\usepackage{color}
\usepackage{amssymb}
\usepackage{amsmath}
\usepackage[english]{babel}
\usepackage{fancyhdr}
\usepackage{amsmath}

\newtheorem{thm}{Theorem}[section]
\newtheorem{lemma}{Lemma}[section]
\newtheorem{prop}{Proposition}[section]
\newtheorem{cor}{Corollary}[section]
\newtheorem{defn}{Definition}[section]
\newtheorem{rem}{Remark}[section]

\def\ep{\hfill $\Diamond$}
\def\R{{\mathfrak R}\, }
\def\S{{\mathfrak S}\, }
\def\M{{\mathfrak M}\, }
\def\T{{\mathfrak T}\, }
\def\Z{{\mathfrak Z}\, }
\def\ci{\begin{color}{red}\,}
\def\cf{\end{color}\,}

\begin{document}
\begin{center}
{\Large\bf Lie $n-$multiplicative mapping on Triangular $n-$Matrix Rings}

\vspace{.2in}
{\bf Bruno L. M. Ferreira}
\\
\textbf{and}
\\
{\bf Henrique Guzzo Jr.}
%$^{a}$}
\vspace{.2in}

%$^{a}$
Universidade Tecnol\'{o}gica Federal do Paran\'{a}, Avenida Professora Laura Pacheco Bastos, 800, 85053-510, Guarapuava, Brazil
\\
and
\\
Universidade de S\~{a}o Paulo, Instituto de
Matem\'{a}tica e Estat\'{i}stica, Rua do Mat\~{a}o,
1010, 05508-090 - S\~{a}o Paulo, Brazil
\vspace{.2in}

brunoferreira@utfpr.edu.br
\\
and
\\
guzzo@ime.usp.br
\vspace{.2in}

\end{center}

{\bf keywords:} Triangular n-matrix rings, additivity, Lie n-multiplicative maps. 

{\bf Mathematics Subject Classification (2010):} 47L35; 16W25

\begin{abstract}
In this paper we extend to triangular $n$-matrix rings and Lie $n$-multiplicative map a result about Lie multiplicative maps on triangular algebras due to Xiaofei Qi and Jinchuan Hou.
\end{abstract}

\section{Introduction}
Let $\mathfrak{R}$ be an associative ring and $\left[x_1,x_2\right] = x_1x_2-x_2x_1$ denote the usual Lie product of $x_1$ and $x_2$.
Let us define the following sequence of polynomials: $p_1(x) = x$ and $p_n(x_1, x_2, \ldots , x_n) = [p_{n-1}(x_1, x_2, \ldots , x_{n-1}) , x_n]$ for all integers $n \geq 2$. Thus, $p_2(x_1, x_2) = [x_1, x_2], \  p_3 (x_1, x_2, x_3) = [[x_1, x_2] , x_3]$, etc. Let $n \geq 2$ be an integer. Assume that $\mathfrak{S}$ is any ring. A map
$\varphi : \R \rightarrow \S$ is called a \textit{Lie $n$-multiplicative mapping} if
\begin{eqnarray}\label{ident1}
\varphi(p_n (x_1, x_2, . . . , x_n)) = p_n (\varphi(x_1), \varphi(x_2), . . . , \varphi(x_n)) 
\end{eqnarray}
In particular, if $n = 2$, $\varphi$ will be called a \textit{Lie multiplicative mapping}. And, if $n=3$, $\varphi$ will be called a \textit{Lie triple multiplicative mapping}.  
\begin{comment}
A mapping $\varphi$ from a ring $\mathfrak{R}$ into another ring $\mathfrak{S}$ is called a Lie triple multiplicative mapping if 
\begin{eqnarray}\label{ident}
\varphi(\left[\left[A,B\right],C\right]) = \left[\left[\varphi(A),\varphi(B)\right],\varphi(C)\right] 
\end{eqnarray}
holds true for all $A,B,C \in \mathfrak{R}$, where $\left[A,B\right] = AB-BA$ is the usual Lie product of $A$ and $B$.
\end{comment}

The study on the question of when a particular application between two rings is additive has become an area of great interest in the theory of rings.
One of the first results ever recorded was given by Martindale III which in his condition requires that the ring possess idempotents, see \cite{Mart}. 
Xiaofei Qi and Jinchuan Hou \cite{QiHou} also considered this question in the context of triangular algebras. They proved the following theorem.

\begin{thm}\cite[Xiaofei Qi and Jinchuan Hou]{QiHou} Let $\mathcal{A}$ and $\mathcal{B}$ be unital algebras over a commutative ring $\mathcal{R}$, and M be a $(\mathcal{A}, \mathcal{B})$-bimodule, which is faithful as a left $\mathcal{A}$-module and also as a right $\mathcal{B}$-module. Let $\mathcal{U} = Tri(\mathcal{A},\mathcal{M}, \mathcal{B})$ be the triangular algebra and $\mathcal{V}$ any algebra over $\mathcal{R}$. Assume that
$\Phi:\mathcal{U} \rightarrow \mathcal{V}$ is a Lie multiplicative isomorphism, that is, $\Phi$ satisfies
$$\Phi(ST - TS) = \Phi(S)\Phi(T) - \Phi(T)\Phi(S) \  \  \   \forall S, T \in \mathcal{U}.$$
Then $\Phi(S + T) = \Phi(S) + \Phi(T) + Z_{S,T}$ for all $S, T \in \mathcal{U}$, where $Z_{S,T}$ is an element in the
centre $\mathcal{Z}(\mathcal{V})$ of $\mathcal{V}$ depending on $S$ and $T.$
\end{thm}

This motivated us to discuss the additivity of Lie $n$-multiplicative mapping on another kind of rings: triangular $n$-matrix rings. In this paper, we give a full answer for this discuss, where the result Xiaofei Qi and Jinchuan Hou is a consequence of our case.
 
\section{Motivation and Definition}

For any unital ring $\R$, let Mod($\R$) denote the category of unitary $\R$-modules, i.e. satisfying $1m = m$ for all elements $m$. This category is important in many areas of mathematics such as ring theory, representation theory and homological algebra.
The purpose of this paper is to work with more general category Mod($\R$) to nonunital rings $\R$. 
It is worth noting that if $\R$ be a nonunital ring, and $\widetilde{\R}$ denote the unital ring $\R \times \mathbb{Z}$ obtained by adjoining an identity. Define operations on $\widetilde{\R}$ by
$$(r, \lambda)+(t, \mu):=(r+t, \lambda + \mu)$$
$$(r, \lambda)\cdot(t, \mu):=(rt + \lambda t + \mu r , \lambda \mu).$$
Then $\widetilde{\R}$ is a ring with $(0_{\R},1) := 1_{\R}$ as multiplicative identity. If $M$ is an non unitary $\R$-module, define an right $\widetilde{\R}$-module operation by 
$$(r, \lambda)m := rm + \lambda m$$
and an left $\widetilde{\R}$-module operation by 
$$m(r, \lambda) := mr + \lambda m,$$
where the action $\mathbb{Z}$ on $M$ is the usual of $M$ as an $\mathbb{Z}$-module.
A module over $\R$ is the same thing as a unitary $\widetilde{\R}$-module.

The following definition is a generalization of the definition that arises in the work of W. S. Cheung \cite{CHEUNG}. This definition appears in Ferreira's paper \cite{bru}.
 
\begin{defn}\label{pri}
Let $\R_1, \R_2, \cdots, \R_n$ be rings and $\M_{ij}$ $(\R_i, \R_j)$-bimodules with $\M_{ii} = \R_i$ for all $1 \leq i \leq j \leq n.$ Let $\varphi_{ijk}: \M_{ij} \otimes_{\R_j} \M_{jk} \longrightarrow \M_{ik}$ be $(\R_i, \R_k)$-bimodules homomorphisms with $\varphi_{iij}: \R_i \otimes_{\R_i} \M_{ij} \longrightarrow \M_{ij}$ and $\varphi_{ijj}: \M_{ij} \otimes_{\R_j} \R_j \longrightarrow \M_{ij}$ the canonical multiplication maps for all $1 \leq i \leq j \leq k \leq n.$ Write $a b = \varphi_{ijk}(a \otimes b)$ for $a \in \M_{ij},$ $b \in \M_{jk}.$ We consider
\begin{enumerate}
\item[{\it (i)}] $\M_{ij}$ is faithful as a left $\R_i$-module and faithful as a right $\R_j$-module $i<j.$
\item[{\it (ii)}] if $m_{ij} \in \M_{ij}$ is such that $\R_i m_{ij} \R_j = 0$ then $m_{ij} = 0$ $i<j.$
\end{enumerate}
Let {\allowdisplaybreaks\begin{eqnarray*}\allowdisplaybreaks \T = \left\{\left(
\begin{array}{cccc}
r_{11} & m_{12} & \ldots & m_{1n}\\
 & r_{22} & \ldots & m_{2n}\\
 &  & \ddots & \vdots\\
 &  &  & r_{nn}\\
\end{array}
\right)_{n \times n} :\underbrace{ r_{ii} \in \R_{i} ~(= \M_{ii}), ~ m_{ij} \in \M_{ij}}_{(1\leq i< j\leq n)}
\right\}\end{eqnarray*}}
be the set of all $n \times n$ matrices $[m_{ij}]$ with the $(i, j)$-entry $m_{ij} \in \M_{ij}$ for all $1 \leq i \leq j \leq n$. Observe that, with the obvious matrix operations of addition and multiplication, $\T$ is a ring iff $a (b c) = (a b) c$ for all $a \in \M_{ik}$, $b \in \M_{kl}$ and $c \in \M_{lj}$ for all $1 \leq i\leq k\leq l\leq j\leq n.$ When $\T$ is a ring, it is called a \textit{triangular $n$-matrix ring}.
\end{defn}
Note that if $n = 2$ we have the triangular matrix ring. As in \cite{bru} we denote by $ \bigoplus^{n}_{i = 1} r_{ii}$ the element
$\left(\begin{array}{cccc}
r_{11} & 0 & \ldots & 0\\
 & r_{22} & \ldots & 0\\
 &  & \ddots & \vdots\\
 &  &  & r_{nn}\\
\end{array}\right)$
in $\T.$

Set $\T_{ij}= \left\{\left(m_{kt}\right):~ m_{kt} = \left\{{ \begin{matrix} m_{ij}, & \textrm{if}~(k,t)=(i,j)\\ 0, & \textrm{if}~(k,t)\neq (i,j)\end{matrix}}, ~i \leq j \right\}.\right.$ Then we can write $ \T = \bigoplus_{1 \leq i \leq j \leq n}\T_{ij}.$ Henceforth the element $a_{ij}$ belongs $\T_{ij}$ and the corresponding elements are in $\R_1, \cdots, \R_n$ or $\M_{ij}.$ By a direct calculation $a_{ij}a_{kl} = 0$ if $j \neq k.$
Also as in \cite{bru} we define natural projections $\pi_{\R_{i}} : \T \longrightarrow \R_{i}$ $(1\leq i\leq n)$ by \\ \centerline{$\left(\begin{array}{cccc}
r_{11} & m_{12} & \ldots & m_{1n}\\
 & r_{22} & \ldots & m_{2n}\\
 &  & \ddots & \vdots\\
 &  &  & r_{nn}\\
\end{array}\right)\longmapsto r_{ii}.$}

\begin{defn}
Let $\R$, $\S$ be rings, we shall say that the Lie $n$-multiplicative mapping $\varphi : \R \rightarrow \S$ is almost additive if there exist $S_{A,B}$ in the centre $\mathcal{Z}(\S)$ of $\S$ depending on $A$ and $B$ such that
$$\varphi(A + B) = \varphi(A) + \varphi(B) + S_{A,B}$$
for all $A, B \in \R$.
\end{defn}

The proposition below appears in \cite{bru} is a generalization of Proposition $3$ of \cite{CHEUNG} and will be very useful.

\begin{prop}\label{seg}
Let $\T$ be a triangular $n-$matrix ring. The center of $\T$ is \\
\centerline{$\mathfrak{Z}(\T) = \left\{ \bigoplus_{i=1}^{n} r_{ii} ~\Big|~ r_{ii}m_{ij} = m_{ij}r_{jj} \mbox{ for all }  m_{ij} \in \M_{ij}, ~i < j\right\}.$}\\
Furthermore, $\mathfrak{\Z}(\T)_{ii} \cong \pi_{\R_i}(\mathfrak{Z}(\T))\subseteq \mathfrak{\Z}(\R_i)$, and there exists a unique ring
isomorphism $\tau^j_{i}$ from $\pi_{\R_i}(\Z(\T))$ to $\pi_{\R_j}(\Z(\T))$ $i \neq j$ such that $r_{ii}m_{ij} = m_{ij}\tau^j_{i}(r_{ii})$ for all $m_{ij} \in \M_{ij}.$
\end{prop}

\begin{rem}\label{obsut}
Throughout this paper we shall make some identifications for example:
Let $r_{kk} \in \R_k$ and $m_{ij} \in \M_{ij}$ then 
$$r_{kk} \equiv \left(\begin{array}{ccccccc}
0 & 0 & \ldots & \ldots & 0 & \ldots & 0\\

 & \ddots & \ldots & \ldots & \vdots& \ldots & \vdots\\
 &  &0  & \ldots & 0 & \ldots & 0 \\
&  &  & r_{kk} & 0& \ldots & 0 \\
&  &  & & 0& \ldots & 0\\
&  &  & & &\ddots& \vdots \\
&  &  & & & & 0 \\
\end{array}\right) $$ \\
and
\\  
$$m_{ij} \equiv \left(\begin{array}{ccccccc}
0 & 0 & \ldots & \ldots & 0 & \ldots & 0\\

 & \ddots & \ldots & \ldots & \vdots& \ldots & \vdots\\
 &  &0  & \ldots & m_{ij} & \ldots & 0 \\
&  &  & \ddots & \vdots&  & \vdots \\
&  &  & & 0& \ldots & 0\\
&  &  & & &\ddots& \vdots \\
&  &  & & & & 0 \\
\end{array}\right)$$
where $i,j,k \in \left\{1, 2 , \ldots, n\right\}.$

In addition we have the following identifications:

Let $\widetilde{\R} = \R \times \mathbb{Z}$ be a unital ring. If $r \in \R$ then $r \equiv (r, 0)$. And $\varphi \times Id : \R \times \mathbb{Z} \rightarrow \S \times \mathbb{Z}$ with $(\varphi \times Id)(r, \lambda) = (\varphi(r), Id(\lambda))= (\varphi(r), \lambda)$ where $Id$ is identity map on $\mathbb{Z}$. For straightforward calculus it is shown that $\varphi \times Id$ is a Lie $n$-multiplicative mapping. In sometimes we shall do $\varphi \times Id \equiv \varphi$.

\end{rem}
%\begin{rem}\label{anelsem}
%Throughout this paper we shall make identification $(0_{\R}, 1) \equiv 1_{\R}$.

\section{A key Lemma}

In this section the following results are generalizations of those that appear in \cite{QiHou}.

\begin{lemma}\label{key}[key Lemma]
Let $\R_1, \R_2, \cdots, \R_n$ be rings and $\M_{ij}$ $(\R_i, \R_j)$-bimodules, $\M_{ij}$ is faithful as a left $\R_i$-module and faithful as a right $\R_j$-module $i<j$. Let $\T$ be the triangular $n$-matrix ring. Assume that $r_{ii} \in \R_{i}$. If $r_{ii}m_{ij} = m_{ij}r_{jj}$ for all $m_{ij} \in \M_{ij}$, $i < j$, then $r_{ii} \in \mathcal{Z}(\R_{i})$. Furthermore, $ \bigoplus^{n}_{i = 1} r_{ii} \in \mathcal{Z}(\T)$, the centre of $\T$.
\end{lemma}
\proof Let $r'_{ii} \in \R_{i}$ and $m_{ij} \in \M_{ij}$. By hypothesis $r_{ii}m_{ij} = m_{ij}r_{jj}$ for all $m_{ij} \in \M_{ij}$, $i < j$ we get
$$(r'_{ii}r_{ii} - r_{ii}r'_{ii})m_{ij} = r'_{ii}(r_{ii}m_{ij}) - r_{ii}(r'_{ii}m_{ij}) = r'_{ii}m_{ij}r_{jj} - r'_{ii}m_{ij}r_{jj} = 0.$$  
Since $\M_{ij}$ is faithful as a left $\R_i$-module we have $r'_{ii}r_{ii} - r_{ii}r'_{ii} = 0$ for all $r'_{ii} \in \R_{i}$. Hence $r_{ii} \in \mathcal{Z}(\R_{i})$ for $i= 1, \ldots, n-1$. In which case that $i = n$ just use the fact that $\M_{ij}$ is faithful as a right $\R_j$-module. Indeed, let $r'_{nn} \in \R_{n}$ and $m_{hn} \in \M_{hn}$. Again by hypothesis $r_{hh}m_{hn} = m_{hn}r_{nn}$ for all $m_{hn} \in \M_{hn}$, $h < n$ we get
\begin{eqnarray*}
m_{hn}(r'_{nn}r_{nn} - r_{nn}r'_{nn}) &=& (m_{hn}r'_{nn})r_{nn} - (m_{hn}r_{nn})r'_{nn} \\&=& r_{hh}m_{hn}r'_{nn} - r_{hh}m_{hn}r'_{nn} = 0.
\end{eqnarray*}
Therefore, $r'_{nn}r_{nn} - r_{nn}r'_{nn}= 0$ for all $r'_{nn} \in \R_{n}$  it follows that $r_{nn} \in \mathcal{Z}(\R_{n})$.
Now using Proposition \ref{seg}, the centre $\mathcal{Z}(\T)$ of $\T$ is
\\
\centerline{$\mathcal{Z}(\T) = \left\{ \bigoplus_{i=1}^{n} r_{ii} ~\Big|~ r_{ii}m_{ij} = m_{ij}r_{jj} \mbox{ for all }  m_{ij} \in \M_{ij}, ~i < j\right\}.$}
Hence $ \bigoplus^{n}_{i = 1} r_{ii} \in \mathcal{Z}(\T)$. \ep

\begin{lemma}\label{stan}[standard Lemma]
Let $A, B, C \in \R$ and $\varphi(C) = \varphi(A) + \varphi(B)$. Then for any $T_1, T_2, \ldots, T_{n-1} \in \R$, we have 
\begin{eqnarray*}
\varphi(p_n(C, T_1, T_2, \ldots, T_{n-1}))&=&\varphi(p_n(A, T_1, T_2, \ldots, T_{n-1}))\\&+&\varphi(p_n(B, T_1, T_2, \ldots, T_{n-1})). 
\end{eqnarray*}
 \end{lemma}
\proof Using the (\ref{ident1}) we have
\begin{eqnarray*}
\varphi(p_n(C, T_1, T_2, \ldots, T_{n-1})) &=& p_n(\varphi(C), \varphi(T_1), \varphi(T_2), \ldots, \varphi(T_{n-1})) \\&=& p_n(\varphi(A) + \varphi(B), \varphi(T_1), \varphi(T_2), \ldots, \varphi(T_{n-1})) \\&=& p_n(\varphi(A), \varphi(T_1), \varphi(T_2), \ldots, \varphi(T_{n-1})) \\&+& p_n(\varphi(B), \varphi(T_1), \varphi(T_2), \ldots, \varphi(T_{n-1}))\\&=& \varphi(p_n(A, T_1, T_2, \ldots, T_{n-1})) \\&+& \varphi(p_n(B, T_1, T_2, \ldots, T_{n-1})).
\end{eqnarray*}
\ep

Note that if 
\begin{eqnarray*}
&&\varphi(p_n(A, T_1, T_2, \ldots, T_{n-1}))+\varphi(p_n(B, T_1, T_2, \ldots, T_{n-1})) = \\&&\varphi(p_n(A, T_1, T_2, \ldots, T_{n-1}) + p_n(B, T_1, T_2, \ldots, T_{n-1})),
\end{eqnarray*}
then by the injectivity of $\varphi$, we get
$$p_n(C, T_1, T_2, \ldots, T_{n-1})=p_n(A, T_1, T_2, \ldots, T_{n-1}) + p_n(B, T_1, T_2, \ldots, T_{n-1}).$$

\begin{lemma}\label{chave}
Let $\R_1, \R_2, \cdots, \R_n$ be rings and $\M_{ij}$ $(\R_i, \R_j)$-bimodules as in Definition \ref{pri}. 
%$\M_{ij}$ is faithful as a left $\R_i$-module and faithful as a right $\R_j$-module $i<j$.
Let $\T$ be the triangular $n$-matrix ring. If $p_2(A, \T) \in \mathcal{Z}(\T)$ then $A \in \mathcal{Z}(\T)$ for each $A \in \T$. Moreover, if $p_{n}(A, \T, \ldots, \T) \in \mathcal{Z}(\T)$ then $A \in \mathcal{Z}(\T)$ for each $A \in \T$.
\end{lemma}
\proof Let $A = \displaystyle \bigoplus_{1 \leq i \leq j \leq n}a_{ij}$ such that $p_2(A, \T) \in \mathcal{Z}(\T)$. Thus $p_2(A, \R_k) \in \mathcal{Z}(\T)$ for $k = 1, \ldots, n-1$, follows that 
$$p_2(A, \R_k) = \bigoplus_{i}p_2(a_{ik}, \R_k) + \bigoplus_{j}p_2(a_{kj}, \R_k) + p_2(a_{kk}, \R_k),$$
because 
\begin{itemize}
	\item If $i<j$, $k\neq i$ and $k \neq j$ then $p_2(a_{ij}, \R_k) = 0;$
	\item If $i = j$ and $k \neq i$ then $p_2(a_{ii}, \R_k) = 0.$
\end{itemize}
As $p_2(A, \R_{k}) \in \mathcal{Z}(\T)$ we have $p_2(A, \R_{k})=p_2(a_{kk}, \R_{k})$ for $k = 1, \ldots , n-1$ by Proposition \ref{seg}.
Consequently, $p_2(A,\R_k)m_{kj} = 0$ for $k < j$. Since $\M_{kj}$ is faithful as a left $\R_k$-module, we see that $p_2(A,\R_k) = 0$ for $k= 1, \ldots , n-1$. 
In the case $k = n$ we use analogous argument and we obtain $m_{in}p_2(A,\R_n) = 0$. And since $\M_{in}$ is faithful as a right $\R_n$-module, it follows that $p_2(A,\R_n) = 0$.
   
Now note that $p_2(A, \M_{ij}) = p_{2}(A_{ii}, \M_{ij}) + p_2(A_{jj}, \M_{ij})$ and as  \linebreak $p_2(A, \M_{ij}) \in \mathcal{Z}(\T)$ 
we get $p_2(A, \M_{ij}) = 0$ for $i < j$.
Therefore $p_2(A, \T) = 0$, that is, $A \in \mathcal{Z}(\T)$. \ep

\begin{comment}
\proof Let $e_i$ the unity of $\R_i$. As $p_2(A,e_i) = e_ip_2(A,e_i)e_j$ and $p_2(A,\T) \in \mathcal{Z}(\T)$ we get $p_2(A,e_i) = 0$. Now as $$p_2(A,m_{ij}) = p_2(A,e_im_{ij}e_j) = e_ip_2(A,m_{ij})e_j,$$ we have $p_2(A,m_{ij}) = 0$ because $p_2(A,m_{ij}) \in \mathcal{Z}(\T)$.

Consequently, $p_2(A,r_i)m_{ij} = p_2(A,r_im_{ij}) - r_ip_2(A,m_{ij}) = 0$ for $i<j$. Since $\M_{ij}$ is faithful as a left $\R_i$-module, we see that $p_2(A,r_i) = 0$. In the case $$m_{in}p_2(A,r_n) = p_2(A,m_{in}r_n) - p_2(A,m_{in})r_n,$$ since $p_2(A,m_{ij}) = 0$ and $\M_{in}$ is faithful as a right $\R_n$-module, we have $p_2(A,r_n) = 0$. Therefore $A \in \mathcal{Z}(\T)$ for each $A \in \T$.
Finally if $$p_{n}(A, \T, \ldots, \T) \in \mathcal{Z}(\T),$$ 
since $p_n(x_1, x_2, \ldots , x_n) = [p_{n-1}(x_1, x_2, \ldots , x_{n-1}) , x_n]$ for all integers $n \geq 2$ follows that $A \in \mathcal{Z}(\T)$. \ep
\end{comment}
\section{Main results}

Let's state our main result in this section which is a generalization of Theorem $2.1$ in \cite{QiHou}.

\begin{thm}\label{priT}
Let $\T$ be the triangular $n$-matrix ring and $\S$ any ring. Consider $\varphi : \T \rightarrow \S$ a bijection Lie $n$-multiplicative mapping satisfying 
%$\varphi(\mathcal{Z}(\T)) \subset \mathcal{Z}(\S)$,
\begin{enumerate}
\item[\it (i)] $\varphi(\mathcal{Z}(\T)) \subset \mathcal{Z}(\S)$
%\item[\it (ii)] $p_{2}(a_i,\R_i) \in \mathcal{Z}(\R_i)$ then $a_i \in \mathcal{Z}(\R_i)$ for each $a_i \in \R_i$ \ $(i = 1, \ldots, n)$
%\item[\it (ii)] $\R_i \M_{ij} = \M_{ij}$ and $\M_{ij}\R_j = \M_{ij}$ \ \ \ for  \ \ $1\leq i < j \leq n$,
\end{enumerate}
then  $$\varphi(A+B) = \varphi(A) + \varphi(B) + S_{A,B}$$ for all $A, B \in \T$, where $S_{A,B}$ is an element in the centre $\mathcal{Z}(\S)$ of $\S$ depending on $A$ and $B$. 
\end{thm}
To prove the Theorem \ref{priT} we introduced a set of lemmas where almost all are generalizations of claims in \cite{QiHou}. We begin with the following lemma

\begin{lemma}\label{lem1}
$\varphi(0) = 0.$
\end{lemma}
\proof
Indeed, $\varphi(0)=\varphi(p_n(0,0, \ldots , 0)) = p_n(\varphi(0),\varphi(0), \ldots , \varphi(0)) \\ =
[\cdots[[\varphi(0),\varphi(0)],\varphi(0)], \cdots ] =0$. \ep

\begin{lemma}\label{lem2}
For any $A \in \T$ and any $Z \in \mathcal{Z}(\T)$, there exists $S \in \mathcal{Z}(\S)$ such that
$\varphi(A + Z) = \varphi(A) + S$.
\end{lemma}
\proof 
Note that $\varphi^{-1}$ is also bijection Lie $n$-multiplicative map. Let $A \in \T$, $Z' \in \mathcal{Z}(\T)$ and $T_2, \ldots , T_n \in \T$. As $\varphi^{-1}$ is surjective we have $\varphi^{-1}(S') = A$ and $\varphi^{-1}(S) = Z'$. 
Now as $\varphi(\mathcal{Z}(\T)) \subset \mathcal{Z}(\S)$ by condition $(i)$ of Theorem \ref{priT} we have,
\begin{eqnarray*}
\varphi(p_n(\varphi^{-1}(S' + S), T_2, \ldots , T_n)) &=& p_n(S' + S, \varphi(T_2), \ldots, \varphi(T_n)) \\&=& p_n(S', \varphi(T_2), \ldots, \varphi(T_n))\\&=& \varphi(p_n(\varphi^{-1}(S'), T_2, \ldots , T_n)). 
\end{eqnarray*}
Soon, $p_n(\varphi^{-1}(S' + S) - \varphi^{-1}(S'), T_2, \ldots , T_n) = 0$ follows that $\varphi^{-1}(S' + S) - \varphi^{-1}(S') \in \mathcal{Z}(\T)$ by Lemma \ref{chave}. Therefore $\varphi(A + Z) = \varphi(A) + S$. 
\ep 
 
\begin{lemma}\label{lem3}
For any $a_{kk} \in \R_k$, $a_{ii} \in \R_i$, $a_{jj} \in \R_j$ and $m_{ij} \in \M_{ij}$, $i<j$, there exist $S, S_1, S_2 \in \mathcal{Z}(\S)$ such that
\begin{itemize}
  \item $ \displaystyle \varphi(a_{kk} + \bigoplus_{i<j} m_{ij}) = \varphi(a_{kk}) + \varphi(\bigoplus_{i<j}m_{ij}) + S$
	\item $\varphi(a_{ii} + m_{ij}) = \varphi(a_{ii}) + \varphi(m_{ij}) + S_1$,
	\item $\varphi(a_{jj} + m_{ij}) = \varphi(a_{jj}) + \varphi(m_{ij}) + S_2$.
\end{itemize}
\end{lemma}
\proof 
We shall only prove the first item because the demonstration of the others are similar.
As $\varphi$ is surjective, there is an element $\displaystyle H = \bigoplus_{1 \leq i \leq j \leq n}h_{ij} \in \T$ such that
\begin{eqnarray*}\label{idlem3}
\varphi(H) = \varphi(a_{kk}) + \varphi(\bigoplus_{i<j}m_{ij})
\end{eqnarray*} 
Let $b_{jj} \in \R_j$, $j \neq k$ and $T_3, \cdots, T_n \in \T$ by Lemma \ref{stan} we have
\begin{eqnarray*}
\varphi(p_n(H, b_{jj}, T_3, \cdots, T_n ))&=& \varphi(p_n(a_{kk}, b_{jj}, T_3, \cdots, T_n )) \\&+& \varphi(p_n(\bigoplus_{i<j} m_{ij}, b_{jj}, T_3, \cdots, T_n )) \\&=& \varphi(p_n(m_{ij}, b_{jj}, T_3, \cdots, T_n )) 
\end{eqnarray*}
It follows that $p_n(H - m_{ij}, b_{jj}, T_3, \cdots, T_n) \in \mathcal{Z}(\T)$ 
and by Lemma \ref{chave} we get $p_2(H - m_{ij}, b_{jj} ) \in \mathcal{Z}(\T)$ thus $p_2(h_{jj}, b_{jj}) \in \mathcal{Z}(\T)$. 
%Therefore 
%$h_{jj} \in \mathcal{Z}(\T)$ and 
%by condition $(ii)$ of the Theorem \ref{priT} we have $h_{jj} \in \mathcal{Z}(\R_j)$.
Therefore $\displaystyle \bigoplus_{i<j} p_2(h_{ij} - m_{ij}, b_{jj}) = 0$ that is $(h_{ij} - m_{ij})b_{jj} = 0$ for all $b_{jj} \in \R_j$ and by condition $(ii)$ of the Definition \ref{pri} we get $h_{ij} = m_{ij}$.
Now consider $b_{kj} \in \M_{kj}$ by standard Lemma \ref{stan} we have 
\begin{eqnarray*}
\varphi(p_n(H, b_{kj}, T_3, \ldots, T_n)) &=& \varphi(p_n(a_{kk}, b_{kj}, T_3, \ldots, T_n)) \\&+&\varphi(p_n(\bigoplus_{i<j}m_{ij}, b_{kj}, T_3, \ldots, T_n)) \\&=& \varphi(p_n(a_{kk}, b_{kj}, T_3, \ldots, T_n)) + \varphi(0) \\&=& \varphi(p_n(a_{kk}, b_{kj}, T_3, \ldots, T_n)). 
\end{eqnarray*}
It follows that $p_n(H - a_{kk} , b_{kj}, T_3, \ldots, T_n) = 0$ and by Lemma \ref{chave} we have $p_2(H-a_{kk}, b_{kj}) \in \mathcal{Z}(\T)$. Thus by Proposition \ref{seg} we get $p_2(H-a_{kk}, b_{kj}) = 0$ which imply that
$(h_{kk} - a_{kk})b_{kj} = b_{kj}h_{jj}$ for all $b_{kj} \in \M_{kj}$.
Therefore by Lemma \ref{key} we obtain $\displaystyle \bigoplus_{l \neq \left\{k, j\right\}} h_{ll} + (h_{kk} - a_{kk} + h_{jj}) \in \mathcal{Z}(\T)$.
And finally by Lemma \ref{lem2} we verified that the Lemma is valid. \ep 

\begin{lemma}\label{lem4}
For any $m_{ij}, s_{ij} \in \M_{ij}$ with $(i < j)$, we have $\varphi(m_{ij} + s_{ij}) = \varphi(m_{ij}) + \varphi(s_{ij})$.
\end{lemma}

\proof Firstly we note that for any $m_{ij}, s_{ij} \in \M_{ij}$, $i<j$, the following identity is valid
$$m_{ij} + s_{ij} = p_n(1_{\R_i} +s_{ij}, 1_{\R_j} +  m_{ij}, 1_{\R_j}, \ldots, 1_{\R_j} ).$$
In deed, due to the Remark \ref{obsut} we get $m_{ij} + s_{ij} = p_2(1_{\R_i} +s_{ij}, 1_{\R_j} +  m_{ij})$. It follows that $m_{ij} + s_{ij} = p_n(1_{\R_i} + s_{ij}, 1_{\R_j} + m_{ij}, 1_{\R_j}, \ldots, 1_{\R_j} )$.
Finally by Lemma \ref{lem3} we have,
\begin{eqnarray*}
\varphi(m_{ij} + s_{ij}) &=& \varphi(p_n(1_{\R_i} + s_{ij},1_{\R_j} + m_{ij} , 1_{\R_j}, \ldots, 1_{\R_j})) \\&=& p_n(\varphi(1_{\R_i} + s_{ij}),\varphi(1_{\R_j} + m_{ij}) , \varphi(1_{\R_j}), \ldots, \varphi(1_{\R_j})) \\&=& p_n(\varphi(1_{\R_i}) + \varphi(s_{ij}) + S_1,\varphi(1_{\R_j}) + \varphi(m_{ij}) + S_2 ,\varphi(1_{\R_j}), \ldots, \varphi(1_{\R_j})) \\&=& p_n(\varphi(1_{\R_i}) + \varphi(s_{ij}),\varphi(1_{\R_j}) + \varphi(m_{ij}) ,\varphi(1_{\R_j}), \ldots, \varphi(1_{\R_j})) \\&=& p_n(\varphi(1_{\R_i}),\varphi(1_{\R_j}),\varphi(1_{\R_j}), \ldots, \varphi(1_{\R_j})) \\&+& p_n(\varphi(1_{\R_i}),\varphi(m_{ij}) ,\varphi(1_{\R_j}), \ldots, \varphi(1_{\R_j})) \\&+& p_n(\varphi(s_{ij}),\varphi(1_{\R_j}),\varphi(1_{\R_j}), \ldots, \varphi(1_{\R_j})) \\&+& p_n(\varphi(s_{ij}),\varphi(m_{ij}) ,\varphi(1_{\R_j}), \ldots, \varphi(1_{\R_j})) \\&=& \varphi(p_n(1_{\R_i}, 1_{\R_j}, 1_{\R_j}, \ldots , 1_{\R_j})) + \varphi(p_n(1_{\R_i}, m_{ij}, 1_{\R_j}, \ldots , 1_{\R_j})) \\&+& \varphi(p_n(s_{ij}, 1_{\R_j}, 1_{\R_j}, \ldots , 1_{\R_j})) + \varphi(p_n(s_{ij}, m_{ij}, 1_{\R_j}, \ldots , 1_{\R_j})) \\&=& \varphi(p_n(1_{\R_i}, m_{ij}, 1_{\R_j}, \ldots , 1_{\R_j})) + \varphi(p_n(s_{ij}, 1_{\R_j}, 1_{\R_j}, \ldots , 1_{\R_j})) \\&=& \varphi(0) + \varphi(m_{ij}) + \varphi(s_{ij}) + \varphi(0) \\&=& \varphi(m_{ij}) + \varphi(s_{ij})
\end{eqnarray*}
\ep 
 
\begin{lemma}\label{lem5}
For any $a_{ii}, b_{ii} \in \R_{i}$, $i = 1, 2, \ldots, n$, there exist $S_i \in \mathcal{Z}(\S)$ such that $\varphi(a_{ii} + b_{ii}) = \varphi(a_{ii}) + \varphi(b_{ii}) + S_i$.
\end{lemma}

\proof
As $\varphi$ is surjective, there is an element $\displaystyle H = \bigoplus_{1 \leq i \leq j \leq n}h_{ij} \in \T$ such that
\begin{eqnarray*}\label{idlem5}
\varphi(H) = \varphi(a_{ii}) + \varphi(b_{ii})
\end{eqnarray*} 
Let $c_{kk} \in \R_k$, $k \neq i$ and $T_3, \cdots, T_n \in \T$ by Lemma \ref{stan} we have
\begin{eqnarray*}
\varphi(p_n(H, c_{kk}, T_3, \cdots, T_n ))&=& \varphi(p_n(a_{ii}, c_{kk}, T_3, \cdots, T_n )) \\&+& \varphi(p_n(b_{ii}, c_{kk}, T_3, \cdots, T_n )) \\&=& \varphi(0) + \varphi(0) = 0 
\end{eqnarray*}
It follows that $p_n(H , c_{kk}, T_3, \cdots, T_n) \in \mathcal{Z}(\T)$ 
and by Lemma \ref{chave} we get $p_2(H, c_{kk} ) \in \mathcal{Z}(\T)$ thus $p_2(h_{kk}, c_{kk}) \in \mathcal{Z}(\T)$. Therefore 
%$h_{kk} \in \mathcal{Z}(\T)$ and
by condition $(i)$ of the Definition \ref{pri} we have $h_{kk} \in \mathcal{Z}(\R_k)$.
Moreover $\displaystyle \bigoplus_{i<j} p_2(h_{ij}, c_{kk}) = 0$ that is $h_{ik}c_{kk} = 0$ for all $c_{kk} \in \R_k$ and by condition $(ii)$ of the Definition \ref{pri} we get $h_{ik} = 0$.
Now consider $c_{ij} \in \M_{ij}$ and $r_{jj} \in \R_j$ by standard Lemma \ref{stan} and Lemma \ref{lem4} we have 
\begin{eqnarray*}
\varphi(p_n(H, c_{ij}, r_{jj}, \ldots, r_{jj})) &=& \varphi(p_n(a_{ii}, c_{ij}, r_{jj}, \ldots, r_{jj})) \\&+&\varphi(p_n(b_{ii}, c_{ij}, r_{jj}, \ldots, r_{jj})) \\&=& \varphi(p_n(a_{ii} + b_{ii}, c_{ij}, r_{jj}, \ldots, r_{jj})). 
\end{eqnarray*}
It follows that $p_n(H - (a_{ii} + b_{ii}) , c_{ij}, r_{jj}, \ldots, r_{jj}) = 0$ and by $(ii)$ of the Definition \ref{pri} we get $(h_{ii} - (a_{ii}+b_{ii}))c_{ij} = c_{ij}h_{jj}$ for all $c_{ij} \in \M_{ij}$. Therefore by Lemma \ref{key} we obtain $\displaystyle \bigoplus_{l \neq \left\{i, j\right\}} h_{ll} + h_{ii} - (a_{ii}+b_{ii}) + h_{jj} \in \mathcal{Z}(\T)$.
And finally by Lemma \ref{lem2} we verified that the Lemma is valid. \ep 

\begin{lemma}\label{lem6}
For any $T \in \T$ with $\displaystyle T = \bigoplus_{1 \leq i\leq j \leq n } T_{ij}$, there exist $S \in \mathcal{Z}(\S)$ such that $\displaystyle \varphi(T) = \bigoplus_{1 \leq i\leq j \leq n } \varphi(T_{ij}) + S.$ 
\end{lemma}

\proof 
As $\varphi$ is surjective, there is an element $\displaystyle H = \bigoplus_{1 \leq i \leq j \leq n}H_{ij} \in \T$ such that
\begin{eqnarray*}\label{idlem6}
\varphi(H) = \bigoplus_{t=1}^{n}\varphi(T_{tt}) + \bigoplus_{1\leq i \leq j \leq n}\varphi(T_{ij})
\end{eqnarray*} 
Let $c_{kk} \in \R_k$, $k = 1, 2, \ldots , n$ by Lemma \ref{stan} we have
\begin{eqnarray*}
\varphi(p_n(H, c_{kk}, c_{kk}, \cdots, c_{kk} ))&=& \varphi(p_n(T_{kk}, c_{kk}, c_{kk}, \cdots, c_{kk} )) \\&+& \sum_{i=1}^{k-1}\varphi(p_n(T_{ik}, c_{kk}, c_{kk}, \cdots, c_{kk} )) \\&+& \sum_{j= k+1}^{n}\varphi(p_n(T_{kj}, c_{kk}, c_{kk}, \cdots, c_{kk} )). 
\end{eqnarray*}
Now let $c_{ll} \in \R_l$ with $l \in \left\{1, \ldots , k-1\right\}$ again by Lemma \ref{stan} we get
\begin{eqnarray*}
&&\varphi(p_n(p_n(H, c_{kk}, c_{kk}, \cdots , c_{kk}), c_{ll}, c_{ll}, \cdots, c_{ll})) = \\&&
 \ \ \ \ \ \ \ \ \ \ \ \ \ \ \ \ \ \ \ \ \ \ \ \ \ \ \ \ \ \ \varphi(p_n(p_n(T_{lk}, c_{kk}, c_{kk}, \cdots , c_{kk}), c_{ll}, c_{ll}, \cdots, c_{ll})).
\end{eqnarray*}
Since $\varphi$ is injective we have
\begin{eqnarray*}
&&p_n(p_n(H, c_{kk}, c_{kk}, \cdots , c_{kk}), c_{ll}, c_{ll}, \cdots, c_{ll}) = \\&& 
\ \ \ \ \ \ \ \ \ \ \ \ \ \ \ \ \ \ \ \ \ \ \ \ \ \ \ \ \ \ p_n(p_n(T_{lk}, c_{kk}, c_{kk}, \cdots , c_{kk}), c_{ll}, c_{ll}, \cdots, c_{ll}).
\end{eqnarray*}
It follows that $p_n(p_n(H - T_{lk}, c_{kk}, c_{kk}, \cdots , c_{kk}), c_{ll}, c_{ll}, \cdots, c_{ll}) = 0$ and by $(ii)$ of the Definition \ref{pri} we obtain $H_{lk} = T_{lk}$.
Again let $c_{qq} \in \R_q$ with $q \in \left\{k+1, \ldots , n\right\}$ by Lemma \ref{stan} we get
\begin{eqnarray*}
&&\varphi(p_n(p_n(H, c_{kk}, c_{kk}, \cdots , c_{kk}), c_{qq}, c_{qq}, \cdots, c_{qq})) = \\&&
 \ \ \ \ \ \ \ \ \ \ \ \ \ \ \ \ \ \ \ \ \ \ \ \ \ \ \ \ \ \ \varphi(p_n(p_n(T_{kq}, c_{kk}, c_{kk}, \cdots , c_{kk}), c_{qq}, c_{qq}, \cdots, c_{qq})).
\end{eqnarray*}
Since $\varphi$ is injective we have
\begin{eqnarray*}
&&p_n(p_n(H, c_{kk}, c_{kk}, \cdots , c_{kk}), c_{qq}, c_{qq}, \cdots, c_{qq}) = \\&& 
\ \ \ \ \ \ \ \ \ \ \ \ \ \ \ \ \ \ \ \ \ \ \ \ \ \ \ \ \ \ p_n(p_n(T_{kq}, c_{kk}, c_{kk}, \cdots , c_{kk}), c_{qq}, c_{qq}, \cdots, c_{qq}).
\end{eqnarray*}
It follows that $p_n(p_n(H - T_{kq}, c_{kk}, c_{kk}, \cdots , c_{kk}), c_{qq}, c_{qq}, \cdots, c_{qq}) = 0$ and by $(ii)$ of the Definition \ref{pri} we obtain $H_{kq} = T_{kq}$.
Finally let $c_{tt} \in \R_{t}$ and $c_{kt} \in \M_{kt}$, $k < t$ by Lemma \ref{stan} we have
\begin{eqnarray*}
&&\varphi(p_n(H, c_{kt}, c_{tt}, \ldots , c_{tt})) = \varphi(p_n(T_{kk}, c_{kt}, c_{tt}, \ldots , c_{tt})) \\&+& \varphi(p_n(T_{tt}, c_{kt}, c_{tt}, \ldots , c_{tt})) + \sum_{i=1}^{k-1} \varphi(p_n(T_{ik}, c_{kt}, c_{tt}, \ldots , c_{tt})) \\&+& \sum_{j=t+1}^{n} \varphi(p_n(T_{tj}, c_{kt}, c_{tt}, \ldots , c_{tt})) = 
\varphi(p_n(T_{kk}, c_{kt}, c_{tt}, \ldots , c_{tt})) \\&+& \varphi(p_n(T_{tt}, c_{kt}, c_{tt}, \ldots , c_{tt})) + \sum_{i=1}^{k-1} \varphi(p_n(T_{ik}, c_{kt}, c_{tt}, \ldots , c_{tt})) + \varphi(0) \\&=& \varphi(p_n(T_{kk}, c_{kt}, c_{tt}, \ldots , c_{tt})) + \varphi(p_n(T_{tt}, c_{kt}, c_{tt}, \ldots , c_{tt})) \\&+& \sum_{i=1}^{k-1} \varphi(p_n(T_{ik}, c_{kt}, c_{tt}, \ldots , c_{tt})).
\end{eqnarray*} 
Now let $c_{kk} \in \R_{k}$ by Lemma \ref{stan} and Lemma \ref{lem4} we obtain
\begin{eqnarray*}
&&\varphi(p_n(p_n(H, c_{kt}, c_{tt}, \ldots , c_{tt}),c_{kk}, c_{kk}, \ldots , c_{kk})) \\&=& \varphi(p_n(p_n(T_{kk}, c_{kt}, c_{tt}, \ldots , c_{tt}),c_{kk}, c_{kk}, \ldots , c_{kk})) \\&+& \varphi(p_n(p_n(T_{tt}, c_{kt}, c_{tt}, \ldots , c_{tt}),c_{kk}, c_{kk}, \ldots , c_{kk})) \\&+& \sum_{i=1}^{k-1} \varphi(p_n(p_n(T_{ik}, c_{kt}, c_{tt}, \ldots , c_{tt}), c_{kk}, c_{kk} , \ldots , c_{kk})) = \\&=& \varphi(p_n(p_n(T_{kk}, c_{kt}, c_{tt}, \ldots , c_{tt}),c_{kk}, c_{kk}, \ldots , c_{kk}) \\&+& p_n(p_n(T_{tt}, c_{kt}, c_{tt}, \ldots , c_{tt}),c_{kk}, c_{kk}, \ldots , c_{kk})) + \sum_{i=1}^{k-1} \varphi(0) \\&=& \varphi(p_n(p_n(T_{kk} + T_{tt}, c_{kt}, c_{tt}, \ldots , c_{tt}),c_{kk}, c_{kk}, \ldots , c_{kk})).
\end{eqnarray*}  
Since $\varphi$ is injective we have
$$p_n(p_n(H - (T_{kk} + T_{tt}), c_{kt}, c_{tt}, \ldots , c_{tt}),c_{kk}, c_{kk}, \ldots , c_{kk}) = 0.$$
By $(ii)$ of the Definition \ref{pri} it follows that $(H_{kk} - T_{kk})c_{kt} = c_{kt}(H_{tt} - T_{tt})$ for all $c_{kt} \in \M_{kt}$.
Therefore $\displaystyle \bigoplus_{i=1}^{n} H_{ii} = \bigoplus_{i=1}^{n} T_{ii}  +  Z$ where $Z \in \mathcal{Z}(\T)$. Now by Lemma \ref{lem2} the result is true. \ep  

\vspace{0,5cm}

We are ready to prove our Theorem \ref{priT}.

\vspace{0,5cm}

\noindent Proof of Theorem. Let $A, B \in \T$. By previous Lemmas we have

\begin{eqnarray*}
\varphi(A + B) &=& \varphi(\bigoplus_{1\leq i \leq j \leq n} A_{ij} + \bigoplus_{1\leq i \leq j \leq n} B_{ij}) \\&=& \varphi(\bigoplus_{k=1}^{n} (A_{kk} + B_{kk}) + \bigoplus_{1\leq i \leq j \leq n} (A_{ij} + B_{ij})) \\&=& \bigoplus_{k=1}^{n}\varphi((A_{kk}+B_{kk})) + \bigoplus_{1\leq i \leq j \leq n}\varphi((A_{ij} + B_{ij})) + S'_{A,B} \\&=& \bigoplus_{k=1}^{n}\varphi(A_{kk}) + \bigoplus_{k=1}^{n}\varphi(B_{kk}) + \bigoplus_{k=1}^{n}S_{\R_{k}} + \bigoplus_{1\leq i \leq j \leq n}\varphi(A_{ij}) \\&+&  \bigoplus_{1\leq i \leq j \leq n}\varphi(B_{ij}) + S'_{A,B} \\&=& \varphi(\bigoplus_{k=1}^{n} A_{kk} + \bigoplus_{1\leq i \leq j \leq n}\varphi(A_{ij})) - S_{A} \\&+& \varphi(\bigoplus_{k=1}^{n} B_{kk} + \bigoplus_{1\leq i \leq j \leq n}\varphi(B_{ij})) - S_{B} + \bigoplus_{k=1}^{n}S_{\R_{k}} + S'_{A,B} \\&=& \varphi(A) + \varphi(B) + S_{A,B},
\end{eqnarray*}
where $\displaystyle S_{A,B} = \bigoplus_{k=1}^{n}S_{\R_{k}} - S_{A} - S_{B} + S'_{A,B}$, so the Theorem \ref{priT} is proved. \ep 

\section{Final Remarks}

\begin{cor}
Let $\T$ be the triangular $n$-matrix unital ring and $\S$ be the ring satisfying
\begin{itemize}
	\item If $p_2(s, \S) \in \mathcal{Z}(\S)$ then $s \in \mathcal{Z}(\S).$ 
\end{itemize}
Then any bijective Lie $n$-multiplicative mapping is almost additive.
\end{cor}
\proof
In deed, let $\varphi: \T \rightarrow \S$ a bijective Lie $n$-multiplicative mapping and $T_2, \ldots, T_n \in \T$ for any $Z \in \mathcal{Z}(\T)$, we have 
\begin{eqnarray*}
p_n(\varphi(Z), \varphi(T_2), \ldots, \varphi(T_n)) = \varphi(p_n(Z, T_2, \ldots, T_n)) = \varphi(0) = 0.
\end{eqnarray*}
Since $T_2, \ldots T_n$ are arbitrary and $\varphi$ is surjective follow that $\varphi(Z) \in \mathcal{Z}(\S)$. \ep 

%The result bellow appear in \cite{posner}
\begin{prop}
For any prime ring $\R$ the following statement
\begin{itemize}
	\item If $p_2(r, \R) \in \mathcal{Z}(\R)$ then $r \in \mathcal{Z}(\R)$
\end{itemize}
is holds true.
\end{prop}
\proof See Lemma $3$ in \cite{posner}. \ep

\begin{thm}
Let $\T$ be the triangular $n$-matrix unital ring and $\S$ be the prime ring. Then any bijective Lie $n$-multiplicative mapping is almost additive.
\end{thm}

\begin{cor}\noindent{ [\textbf{Xiaofei Qia and Jinchuan Hou} \cite{QiHou}]}
Let $\mathcal{A}$ and $\mathcal{B}$ be unital algebras over a commutative ring $\mathcal{R}$, and M be a $(\mathcal{A}, \mathcal{B})$-bimodule, which is faithful as a left $\mathcal{A}$-module and also as a right $\mathcal{B}$-module. Let $\mathcal{U} = Tri(\mathcal{A},\mathcal{M}, \mathcal{B})$ be the triangular algebra and $\mathcal{V}$ any algebra over $\mathcal{R}$. Assume that
$\Phi:\mathcal{U} \rightarrow \mathcal{V}$ is a Lie multiplicative isomorphism, that is, $\Phi$ satisfies
$$\Phi(ST - TS) = \Phi(S)\Phi(T) - \Phi(T)\Phi(S) \  \  \   \forall S, T \in \mathcal{U}.$$
Then $\Phi(S + T) = \Phi(S) + \Phi(T) + Z_{S,T}$ for all $S, T \in \mathcal{U}$, where $Z_{S,T}$ is an element in the
centre $\mathcal{Z}(\mathcal{V})$ of $\mathcal{V}$ depending on $S$ and $T.$
\end{cor}

\proof This is consequence of our Theorem \ref{priT} for $n = 2$. \ep

\section{Application in Nest Algebras}

A \textit{nest} $\mathcal{N}$ is a totally ordered set of closed subspaces of a Hilbert space $\mathcal{H}$ such that $\left\{0\right\}, \mathcal{H} \in \mathcal{N}$ , and $\mathcal{N}$ is closed under the taking of arbitrary intersections and closed linear spans of its elements. The \textit{nest algebra} associated to $\mathcal{N}$ is the set $\mathcal{T}(\mathcal{N}) = \left\{T \in \mathcal{B}(\mathcal{H}) ~:~ TN \subseteq N \mbox{ for all } N \in\mathcal{N}\right\},$ where $\mathcal{B}(\mathcal{H})$ is the algebra of bounded operators over a complex Hilbert space $\mathcal{H}$.

We recall the standard result (\cite{CHEUNG}, Proposition $16$) that say we can view $\mathcal{T}(\mathcal{N})$ as triangular algebra
$\left(\begin{matrix}
A&M\\
 &B
\end{matrix}\right)$
where $A, B$ are themselves nest algebras.
\begin{prop}
If $N \in \mathcal{N} \setminus \left\{0, H\right\}$ and $E$ is the orthonormal projection onto $N$, then $E\mathcal{N}E$ and $(1-E)\mathcal{N}(1-E)$ are nest, $\mathcal{T}(E\mathcal{N}E) = E \mathcal{T}(\mathcal{N})E$ and 
$\mathcal{T}((1-E)\mathcal{N}(1-E)) = (1-E) \mathcal{T}(\mathcal{N})(1-E)$. Furthermore
$$\mathcal{T}(\mathcal{N}) = \left(\begin{matrix}
\mathcal{T}(E\mathcal{N}E)&E\mathcal{T}(\mathcal{N})(1-E)\\
 &\mathcal{T}((1-E)\mathcal{N}(1-E))
\end{matrix}\right). $$
\end{prop}
We refer the
reader to \cite{Davi} for the general theory of nest algebras.

\begin{cor}
Let $P_n$ be an increasing sequence of finite dimensional subspaces such that their union is dense in $\mathcal{H}$. 
Consider $\mathcal{P} = \left\{\left\{0\right\}, P_n, n\geq1, \mathcal{H}\right\}$ a nest and $\mathcal{T}(\mathcal{P})$ the set consists of all operators which have a block upper triangular matrix with respect to $\mathcal{P}$.
If a mapping $\varphi : \mathcal{T}(\mathcal{P}) \longrightarrow \mathcal{T}(\mathcal{P})$ satisfies
$$\varphi(p_2([f,g]) = p_2(\varphi(f),\varphi(g))$$
for all $f, g \in \mathcal{T}(\mathcal{P})$, then $\varphi$ is almost additive.
\end{cor}


\begin{thebibliography}{99}
\bibitem{CHEUNG} W. S. Cheung, {\it Commuting maps of triangular algebras}, J. London Math. Soc. {\bf 63}, 17–127 (2001).
%\bibitem{Daif} M. Daif, {\it When is a multiplicative derivation additive ?}, Internat. J Math. and Math. Sci. {\bf 14}, 615-618 (1991).
\bibitem{Davi} K.R. Davidson, {\it Nest algebras}, Pitman Research Notes in Mathematics Series {\bf 191}, (1988).
\bibitem{bru} B.L.M. Ferreira, {\it Multiplicative maps on triangular n-matrix rings}, International Journal of Mathematics, Game Theory and Algebra, {\bf 23}, 1-14 (2014).
\bibitem{Mart} W. S. Martindale III, {\it When are multiplicative mappings additive?}, Proc. Amer. Math. Soc. {\bf 21}, 695-698 (1969).
\bibitem{posner} E.C. Posner, {\it Derivations in prime rings}, Proc. Amer. Math. Soc. {\bf 8}, 1093-1100 (1957).
\bibitem{QiHou} Xiaofei Qi and Jinchuan Hou, {\it Additivity of Lie multiplicative maps on triangular algebras}, Linear and Multilinear Algebra, {\bf 59}, 391-397 (2011).
%\bibitem{Wangc} Yu Wang, {\it The Additivity of multiplicative maps on rings}, Communications in Algebra {\bf 37} (2009), 2351-2356.
%\bibitem{Wang} Yu Wang, {\it Additivity of multiplicative maps on triangular rings}, Linear Algebra and its Applications {\bf 434}, 625-635 (2011).



\end{thebibliography}
\end{document}